\documentclass[11pt]{article}
\usepackage{amsfonts}
\usepackage{subfigure} 
\usepackage{epsfig}
\usepackage{amsmath, amsthm}
\newcommand{\be}{\begin{eqnarray}}     	\newcommand{\ee}{\end{eqnarray}}

\newcommand{\vol}{\mathrm{Vol}}
\newcommand{\dist}{\mathrm{dist}}
\newcommand{\rem}{\mathrm{Rm}}
\newcommand{\ric}{\mathrm{Ric}}
\newcommand{\diam}{\mathrm{diam}}

\title{Limiting behaviour of the Ricci flow}
\author{Natasa Sesum}
\date{}
\theoremstyle{plain}
\newtheorem{dummy}{Dummy}

\theoremstyle{definition}
\newtheorem{definition}[dummy]{Definition}
\newtheorem{step}{Step}[dummy]

\theoremstyle{plain}

\newtheorem{remark}[dummy]{Remark}
\newtheorem{lemma}[dummy]{Lemma}
\newtheorem{theorem}[dummy]{Theorem}

\newtheorem{claim}[dummy]{Claim}

\begin{document}

\maketitle

\begin{abstract}
We will consider a {\it $\tau$-flow}, given by the equation
$\frac{d}{dt}g_{ij} = -2R_{ij} + \frac{1}{\tau}g_{ij}$ on a closed
manifold $M$, for all times $t\in [0,\infty)$. We will prove that if
the curvature operator and the diameter of $(M,g(t))$ are uniformly
bounded along the flow, then we have a sequential convergence of the
flow toward the solitons.
\end{abstract}

\begin{section}{Introduction}

The studies of singularities and the limiting behaviours of solutions
of various geometric partial differential equations have been
important in geometric analysis. One of these important geometric
equations is so called Ricci flow equation, itroduced by Richard
Hamilton in \cite{hamilton1982}. It is the equation
$\frac{d}{dt}g_{ij}(t) = -2R_{ij}$, for a Riemannian metric
$g_{ij}(t)$. The short time existence of this equation was proved by
Hamilton in \cite{hamilton1982} and somewhat later the proof was
significantly simplified by DeTurck in \cite{deturck1983}. Hamilton
showed that the Ricci flow preserves the positivity of the Ricci
tensor in dimension three and of the curvature operator in all
dimensions. This observation helped him to prove the convergence
results in dimensions three and four, towards metrics of constant
positive curvatures (in the case of positive Ricci curvature and
positive curvature operator respectively).

Besides the short time existence we can also study a long time
existence of the Ricci flow. There is a well known Hamilton's
result.

\begin{theorem}[Hamilton]
For any smooth initial metric on a compact manifold there exists a
maximal time $T$ on which there is a unique smooth solution to the
ricci flow for $0 \le t < T$. Either $T = \infty$ or else the
curvature is unbounded as $t\to T$.
\end{theorem}

One can ask what happens to a solution if it exists for all times and
under which conditions it will converge to a metric that will have
nice properties. In the case of dimension three with positive Ricci
curvature and dimension four with positive curvature operator we know
that a solution converges to an Einstein metric. In general, we can
not expect to get an Einstein metric in the limit. We can expect to
get a solution to an evolution equation which moves under a
one-parameter subgroup of the symmetry group of the equation. These
kinds of solutions are called solitons.

Our goal in this paper is to prove the following theorem.

\begin{theorem}[Main Theorem]
\label{theorem-theorem_main}
Consider the flow
\be
\label{equation-equation_our_flow}
\frac{d g_{ij}}{d t} = -2R_{ij} + \frac{1}{\tau}g_{ij} \ee 
on a
compact manifold $M$, where $\tau > 0$ is fixed, $|Rm| \leq C$ and
$diam(M,g(t) \leq C$ $\forall t\in [0,\infty)$. Then for every
sequence of times $t_i\to\infty$ there exists a subsequence, 
so that $g(t_i+t)\to h(t)$ and $h(t)$ is a Ricci soliton.
\end{theorem}

The organization of the paper is as follows. In section $3$ we will
prove some properties of $\mu(g,\tau)$ that has been introduced by
Perelman in \cite{perelman2002}. They will be useful in the later
sections of the paper. In section $3$ we will prove Theorem
\ref{theorem-theorem_main}.

{\bf Acknowledgements}: The author would like to thank her advisor Gang
Tian for bringing this problem to her attention and for constant help
and support.

\end{section}

\begin{section}{Preliminaries}

Perelman's functional $\mathcal{W}$ and its properties will play an
important role in the proof of Theorem \ref{theorem-theorem_main}. $M$
will always denote a compact manifold, and $(g_{ij})_t = -2R_{ij} +
\frac{1}{\tau}g_{ij}$ will be a flow that we will be considering
throughout the whole paper. Perelman's functional $\mathcal{W}$ has
been introduced in \cite{perelman2002}.

$$W(g,f,\tau) = (4\pi\tau)^{-\frac{n}{2}}\int_M e^{-f}[\tau(|\nabla
f|^2 + R) + f - n] dV_g.$$

We will consider this functional restricted to $f$ satisfying

\begin{equation}
\label{equation-equation_constraint}
\int_M (4\pi\tau)^{-\frac{n}{2}}e^{-f} dV = 1.
\end{equation} 
$\mathcal{W}$ is invariant under simultaneous scaling 
of $\tau$ and $g$.  Perelman showed that the Ricci flow 
can be viewed as a gradient flow of functional $\mathcal{W}$.
Let $\mu(g,\tau) = \inf\mathcal{W}(g,f,\tau)$ over smooth $f$ 
satisfying (\ref{equation-equation_constraint}). It has been 
showed by Perelman that there always exists a smooth minimizer
on a closed manifold $M$, that $\mu(g,\tau)$ is negative 
for small $\tau > 0$ and that it tends to zero as $\tau\to 0$. 
One of the most important properties of $\mathcal{W}$ is
the monotonicity formula.

\begin{theorem}[Perelman]
$\frac{d}{dt}\mathcal{W} = \int_M 2\tau|R_{ij} + \nabla_i\nabla_j f -
\frac{1}{2\tau}g_{ij}|^2(4\pi\tau)^{-\frac{n}{2}}e^{-f} dV \ge 0$ and
therefore $\mathcal{W}$ is increasing along the Ricci flow.
\end{theorem}

One of the very important applications of the monotonicity formula is
noncollapsing theorem for the Ricci flow that has been proved by
Perelman in \cite{perelman2002}.

\begin{definition}
Let $g_{ij}(t)$ be a smooth solution to the Ricci flow
$(g_{ij})_t=-2R_{ij}(t)$ on $[0,T)$. We say that $g_{ij}(t)$ is
loacally collapsing at $T$, if there is a sequence of times $t_k\to T$
and a sequence of metric balls $B_k = B(p_k,r_k)$ at times $t_k$, such
that $\frac{r_k2}{t_k}$ is bounded, $|\rem|(g_{ij}(t_k))\le r_k^{-2}$ in
$B_k$ and $r_k^{-n}\vol(B_k)\to 0$.
\end{definition}

\begin{theorem}
\label{theorem-perelman_theorem}
If $M$ is closed and $T < \infty$, then $g_{ij}(t)$ is not locally
collapsing at $T$.
\end{theorem}

\end{section}

\begin{section}{Sequential convergence of a $\tau$-flow}

\begin{definition}
{\it $\tau$-flow} is given by the equation
\begin{equation}
\frac{d}{dt}g_{ij} = -2R_{ij} + \frac{1}{\tau}g_{ij},
\end{equation}
for $\tau > 0$.
\end{definition}

We want to prove the Theorem \ref{theorem-theorem_main} in this
section.

\begin{subsection}{Convergence toward the solutions of the Ricci flow}

In order to prove Theorem \ref{theorem-theorem_main} we will first
show that it is reasonable to expect a convergence toward a smooth
manifold, i.e. that a limit manifold will not collapse. 

\begin{claim}
\label{claim-claim_noncollapsing}
Consider the flow as above. For every fixed $\tau > 0$ there exists a
constant $C$ such that $\vol_{g(t)}(M)\ge C$ for every $t$, i.e. we
have a uniform lower bound on the volumes.
\end{claim}

\begin{proof}
Assume that the claim is not true, i.e. that there exists a sequence
$t_i$ s.t. $\vol_{g(t_i)}(M)\to 0$ as $i\to \infty$. Let $\bar{g}(s) =
c(s)g(t(s))$ be unnormalized flow, for $s\in [0,\tau)$, where:
$$t(s) = -\tau\ln (1-\frac{s}{\tau}).$$
$$c(s) = 1-\frac{s}{\tau}.$$
$$R(\bar{g}) = \frac{R(g)}{c(s)}.$$
Find $s_i$, such that $t(s_i) = t_i$. We get that $s_i =
\tau(1-e^{-\frac{t_i}{\tau}})$. $s_i\to \tau$ as $i\to \infty$. Let
\be
\max_{M\times [0,s_i]}|Rm|(\bar{g}(s)) = Q_i,
\ee
and assume that the maximum is achieved at $p_i$. By the corollary of
Perelman's noncollapsing theorem we have that:
$$\frac{\vol_{\bar{g}(t)}B(p_i,r)}{r^n}\ge C_1,$$
for $r\le C\sqrt{\frac{\tau}{Q_i}}$ and $t\in [0,s_i)$. Choose $r =
C\sqrt{\frac{\tau}{Q_i}}$ and $t = s_i$. 

$$(\sqrt{Q_i})^n\vol_{\bar{g}(s_i)}B(p_i, C\sqrt{\frac{\tau}{Q_i}})\ge
(C\sqrt{\tau})^n C_1 = \tilde{C}.$$ 
Since $\vol_{\bar{g}(s_i)}B(p_i,r)
= c(s_i)^{\frac{n}{2}}\vol_{g(t_i)}B(p_i,\tilde{r})$, where
$\tilde{r}$ might be a different radius as a matter of scaling and
since $Q_i \le \frac{C}{c(s_i)}$ (because the curvature of $g(t)$ is
uniformly bounded), we get that:
$$\vol_{g(t_i)}(M)\ge \tilde{C}/C,$$
where $\tilde{C}$ and $C$ do not depend on $i$. Let $i\to \infty$ in the
previous inequality to get a contradiction. Therefore we have a
uniform lower bound on volumes.
\end{proof}

\begin{remark}
The assumptions of the Theorem \ref{theorem-theorem_main} and the
result of Claim \ref{claim-claim_noncollapsing} imply the uniform
bounds on the curvature tensors, uniform upper bound on the diameters
and uniform lower bounds on the volumes.  Similarly like in the case
of unnormalized flow, uniform bounds on the curvatures gives us
uniform bounds on all covariant derivatives, so by Hamilton's
compactness theorem, for every sequence $t_i \nearrow \infty$ as $i
\to \infty$, there exists a subsequence (call it again $t_i$), such
that $(M,g(t_i+t))$ converges to $(M,h(t))$, in the sense that there
exist diffeomorphisms $\phi_i : M \to M$, so that $\phi_i ^*g(t_i+t)$
converge uniformly together with their covariant derivatives to
metrics $h(t)$ on compact subsets of $M\times[0,\infty)$. Moreover,
$h(t)$ is a solution of a $\tau$-flow as well.
\end{remark}

\end{subsection}

\begin{subsection}{Continuity of the minimizers for $\mathcal{W}$}

We will recall a definition of Perelman's functional $\mathcal{W} =
(4\pi\tau)^{-\frac{n}{2}}\int_M e^{-f}[\tau(R + |\nabla f|^2) + f -
n]dV$. The constraint on $f$ for this functional is (*)
$(4\pi\tau)^{-\frac{n}{2}}\int e^{-f} dV = 1$. Let $\mu(g,\tau) = \inf
\mathcal{W}(g,f,\tau)$ under the constraint (*).  This infinimum has
been achieved by some smooth minimizer $f$. Perelman has also proved
that for a fixed metric $g$, $\lim_{\tau\to 0}\mu(g,\tau) = 0$ and
$\mu(g,\tau) < 0$ for a small value of $\tau > 0$.

In the case of a $\tau$-flow $g(t)$, $\tau > 0$ is being fixed in
time, and by the monotonicity formula for $\mathcal{W}$ we have that
$\mu(g(t),\tau)$ is increasing along the flow. Therefore, there exists
$\lim_{t\to \infty}\mu(g(t),\tau)$.

\begin{claim} 
\label{claim-step0}
$\lim_{t\to \infty}\mu(g(t),\tau)$ is finite.
\end{claim}

\begin{proof}
Assume that $\lim_{t\to\infty}\mu(g(t),\tau) = \infty$. Then, $\forall
i$, $\exists t_i$ s.t. $\mu(g(t_i),\tau) \ge i$. There exists a
subsequence (call it $t_i$) such that $(M,g_i)$ converges to $(M,h)$,
for some metric $h$. From the first part of Lemma
\ref{lemma-metric_continuity} we get that $\mu(g(t_i),\tau) <
\mu(h,\tau) + \epsilon$, for $i$ big enough. Letting $i\to \infty$ we
get a contradiction.
\end{proof}

\begin{lemma}
\label{lemma-metric_continuity} 
If $(M,g_i)$ tend to $(M,h)$ when $i \to \infty$, where $g_i = g(t_i)$
and $t_i \nearrow \infty$, then $\lim_{i\to \infty}\mu(g_i,\tau) =
\mu(h,\tau)$.
\end{lemma}

\begin{proof}
$$\mu(h,\tau) = \int_M(\tau(|\nabla f |^2 + R(h)) + f -
n)(4\pi\tau)^{-\frac{n}{2}}dV_h.$$ Since $\phi_i^*g_i \to \infty$
uniformly with their covariant derivatives, if $\epsilon >0$ is fixed,
there exists some big $i_0$, so that for $i\ge i_0$
$$\mu(h,\tau) \ge \int_M(\tau(|\nabla f|^2 + R(\tilde{g_i})) + f -
n)(4\pi\tau)^{-\frac{n}{2}}dV_{\tilde{g_i}} - \frac{\epsilon}{2},$$
where $\tilde{g_i} = \phi^*g_i$.
Change the variables in the above integral by diffeomorphism $\phi_i$. 
$$\mu(h,\tau) \ge \int_M(\tau(|\nabla _i f_i|^2 + R(g_i)) + f_i
-n)(4\pi\tau)^{-\frac{n}{2}}dV_{g_i} - \frac{\epsilon}{2},$$ where $f_i
= \phi^* f$.  Perturb a little bit $f_i$ to get $\tilde{f_i}$, by a
quantity that tends to zero, so that $\int_M
e^{-\tilde{f_i}}(4\pi\tau)^{-\frac{n}{2}}dV_{g_i} = 1$. Since our
geometries are uniformly bounded, for big enough $i_0$ we will have
\begin{equation}
\label{equation-inequality1}
\mu(h,\tau) \ge \mathcal{W}(g_i,\tilde{f_i},\tau) - \epsilon \ge
\mu(g_i,\tau) - \epsilon.
\end{equation}
Let $u_i = e^{-\frac{f_i}{2}}$. We have seen that minimizing
$\mu(g_i,\tau)$ by $f_i$ is equivalent to minimizing the
following expression in $u_i$:
$$\int_M\tau(4|\nabla_iu_i|^2 + R_iu_i^2) - 2u_i^2\ln u_i -
nu_i^2)(4\pi\tau)^{-\frac{n}{2}}dV_{g_i}.$$
The minimizer $u_i$ has to satisfy the following elliptic differential equation
\begin{equation}
\label{equation-minimizing_in_u}
\tau(-4\Delta_iu_i + R_iu_i) - 2u_i\ln u_i - nu_i =
\mu_{i,\tau}u_i.
\end{equation}
$\mu_{i,\tau}$ is uniformly bounded, since there is
a finite $\lim_{t\to \infty}\mu(g(t),\tau)$.
Now we can easily get:

\begin{equation}
\label{equation-uniform_bound1}
\int_M u_i^2 (4\pi\tau)^{-\frac{n}{2}}dV_i \le C,
\end{equation}

\begin{equation}
\label{equation-uniform_bound2}
\tau\int_M |\nabla_iu_i|^2(4\pi\tau)^{-\frac{n}{2}}dV_i \le C,
\end{equation}
i.e. $u_i \in W^{1,2}$ with $$||u_i||_{W^{1,2}} \le C \ \ \forall i.$$
From (\ref{equation-minimizing_in_u}), by standard regularity theory of partial
differential equations and Sobolev embedding theorems, we get that
$u_i \in W^{k,p}$ with uniformly bounded $W^{k,p}$ norms, where $p <
\frac{2n}{n-2}$, and therefore with uniformly bounded $C^{2,\alpha}$
norms, i.e. $||u_i||_{C^{2,\alpha}} \le C$.  Furthermore,
\be
\label{equation-equality_help}
\mu(g_i,\tau) = \int_M(\tau(4|\nabla _i u_i|^2 + R_iu_i^2) - 2u_i^2\ln
u_i -nu_i^2)(4\pi\tau)^{-\frac{n}{2}}dV_i\nonumber\\
= \int_M \tau(|\tilde{\nabla}\tilde{u_i}|^2 4 +
\tilde{R_i}\tilde{u_i}^2) - 2\tilde{u_i}\ln \tilde{u_i} -
n\tilde{u_i}^2)(4\pi\tau)^{-\frac{n}{2}}dV_{\tilde{g_i}}, \ee
where $\tilde{u_i} = \phi_i^*u_i$. $\phi_i^*g_i$ is close to $h$ and
therefore for $i$ big enough, $\phi_i$ is almost an isometry, so
$D_j\phi_i^{-1}$ can be uniformly bounded in terms of bounds on $g_i$
and $h$, $g_i$ can be bounded in terms of $h$. We cover $M$ with
finitely many geodesic balls of fixed radius $\rho$ ( we can do it
since we have a uniform bound on the injectivity radii from below). We
use local coordinates in each of the balls to get:
$$ |\tilde{\nabla}_i\tilde{u}_i|^2 = \tilde{g}_i^{jk}D_j(u_i\circ
\phi_i^{-1})D_k(u_i\circ\phi_i^{-1}).$$
$$|\tilde{\nabla}\tilde{u}_i|^2 =
\tilde{g}_i^{jk}(D_ju_i)(D_ku_i)(\phi_i^{-1})D_j\phi_i^{-1}D_k\phi_i^{-1}.$$
Now we can easily conclude that we have a uniform bound on
$|\tilde{\nabla}{u_i}|^2$. Since the integrand in
(\ref{equation-equality_help}) is uniformly bounded in $i$, and since
$\tilde{g}_i$ uniformly converge with their covariant derivatives to
$h$, we have that for $i$ large enough
$$\mu(g_i,\tau) \ge \int_M (\tau(4|{\nabla _h}\tilde{u}_i|^2 +
R_h\tilde{u}_i^2) - 2\tilde{u}_i\ln \tilde{u}_i
-n\tilde{u}_i^2)(4\pi\tau)^{-\frac{n}{2}}dV_h - \epsilon.$$ 
Since $l_i = \int_M \tilde{u}_i^2(4\pi\tau)^{-\frac{n}{2}}dV_h$ is close to $1$
when $i\to \infty$, taking $\bar{u}_i = \frac{\tilde{u}_i}{l_i}$ and
using all the uniform bounds that we have got by now
$$\mu(g_i,\tau)\ge \mathcal{W}(h,\bar{u}_i,\tau) - \epsilon \ge
\mu(h,\tau) - \epsilon.$$ 
By the previous inequality (for $i$ big
enough) and by (\ref{equation-inequality1}) we get $\lim_{i\to
\infty}\mu(g_i,\tau) = \mu(h,\tau)$.
\end{proof}
 
Following the notation from the previous lemma, by Arzela-Ascoli
theorem there exists a subsequence, $u_i$, so that it converges in
$C^{2,\alpha}$ norm to some function $u$. We can also get the higher
order uniform estimates on $u_i$ in a similar manner as in Lemma
\ref{lemma-metric_continuity}. Therefore, to show that a sequence of
minimizers for $\mu(g_i,\tau)$ converges to a minimizer of
$\mu(h,\tau)$ it is enough to show the following lemma.

\begin{lemma}
\label{lemma-lower_bound}
$\exists C>0$ so that $u_i \ge C>0 \ \ \forall i$ and $\forall x\in M$
\end{lemma}

\begin{proof}
Assume that there exists a sequence $u_i$ and $p_i \in M$, such that
$0 < u_i(p_i) <\frac{1}{2i}$. $M$ is compact and therefore there is a
subsequence, $\{p_i\}$ converging to $ p\in M$ when $i \to
\infty$. $C^{2,\alpha}$ norms of $u_i$ are uniformly bounded in $i$
and therefore $u_i(p) < u_i(p_i) + C\dist_{g_i}(p,p_i) \to 0$ as
$i\to\infty$. Let $u$ be a limit of $\{u_i\}$ in $C^{2,\alpha}$
norm. Then $u(p)=0$. Take a geodesic ball $B(p,r)$. Let $f\in
C_0^{\infty}(M)$ be a $C^{\infty}$ function of $r$ alone, compactly
supported in $B(p,r)\backslash \{p\}$.
$$\int_M(\tau(\nabla u_i \nabla f + R_iu_if) -2u_if\ln u_i -nu_if -
\mu(g_i,\tau)u_if)dV_i = 0.$$ 
For this $f$, letting $i\to\infty$, using
the result of the previous lemma and the fact that the integrand in
the previous integral is uniformly bounded in $i$ we get
$$\int_M(\tau(\nabla u\nabla f + fuR(h)) - 2uf\ln u - nuf -
\mu(h,\tau)fu)dV_h = 0.$$
Proceeding in the same manner as in \cite{rothaus1981} we can get that
$u\equiv 0$ in some small ball around $p$. Using the connectedness
argument, $u\equiv 0$ in $M$. On the other hand $\int_M
u_i^2(4\pi\tau)^{-\frac{n}{2}}dV_i = 1$ and letting $i\to \infty$ we
get a contradiction.
\end{proof}

If we write down the equations (\ref{equation-minimizing_in_u}) for
all $\{u_i\}$, letting $i\to \infty$, keeping in mind the previous
lemma we get
$$ \tau(-4\Delta u + R(h))u - 2u\ln u - nu = \mu(h,\tau)u,$$
i.e. $u$ is the minimizer for $\mu(h,\tau)$.

So far we have proved the following theorem
\begin{theorem}
\label{theorem-theorem_continuity_of_minimizers}
If $(M,g_i)\to (M,h)$ as $ i\to \infty$, then for a given $\tau > 0$,
if $\mu(g_i, \tau) = \mathcal{W}(g_i,f_i,\tau)$, then $f_i \to f$ in
$C^{2,\alpha}$ norm, where $\mu(h,\tau) = \mathcal{W}(h,f,\tau)$.
\end{theorem}

\end{subsection}

\begin{subsection}{Further estimates on the minimizers}

In this subsection we want to use the minimizers $f_t$ for $\mathcal{W}$
at different times to construct the functions $f_t(s)$ for $s\in
[0,t]$. By using the parabolic regularity we will be able to get the
uniform estimates on $C^{k,\alpha}$ norms of $f_t(s)$. This will
enable us to take a limit of this functions along the sequences. This
limits are the functions that will turn out to be the potential
functions that come into the equations describing the soliton type
solutions arising in a limit.

For any $t$ we can find
$f_t$ such that $\mathcal{W}(g(t),f_t,\tau) = \mu(g(t),\tau)$. If we
flow $f_t$ backward, we will get functions $f_t(s)$ that satisfy
\begin{eqnarray*}
\frac{df_t(s)}{ds} &=& -R(s) - \Delta f_t(s) + |\nabla f_t(s)|^2 +
\frac{n}{2\tau}, \\
f_t(t) &=& f_t.
\end{eqnarray*}
We know that minimizing $\mathcal{W}$ in $f$ is equivalent to minimizing 
the corresponding functional in
$\tilde{u}$, where $\tilde{u}_t = e^{-\frac{f_t}{2}}$. Let $u_t(s) =
\tilde{u}_t^2(s)$. The equation for $u_t(s)$ is
\begin{eqnarray*}
\frac{du_t}{ds} &=& -\Delta u_t + (-\frac{n}{2\tau} +
R(s))u_t(s), \\
u_t(t) &=& u_t.
\end{eqnarray*}
By the monotonicity of $\mathcal{W}$ along the flow
(\ref{equation-equation_our_flow}) we have that
$$\mu (g(s),\tau ) \le \mathcal{W}(g(s),f_t(s),\tau) \le
\mathcal{W}(g(t),f_t,\tau) = \mu(g(t),\tau).$$ 
First of all, there exists $\lim_{t\to\infty}\mu(g(t),\tau)$. 
It is finite, since for
every sequence $t_i\to\infty$ there exists a subsequence such that
$g(t_i) \to h(0)$ and by Lemma \ref{lemma-metric_continuity} from the
previous section, we have that $\mu(g(t_i),\tau)\to \mu(h(0),\tau)$.

Instead of functional $\mathcal{W}(g(s),f_t(s),\tau)$ we can consider
the equivalent functional which depends on $\tilde{u}_t(s) = e^{-f_t(s)/2}$.
\begin{equation}
\mathcal{W}(u_t(s)) = \int_M[\tau(4|\nabla \tilde{u}_t(s)|^2+R
\tilde{u}_t(s)^2)-\tilde{u}_t(s)^2\log
\tilde{u}_t(s)^2-n\tilde{u}_t(s)^2](4\pi\tau)^{-n/2}dV,
\end{equation}
where $\tilde{u}_t$ satisfy 
$$\tau(-4\Delta \tilde{u}_t + R\tilde{u}_t) -
2\tilde{u}_t\ln\tilde{u}_t - n\tilde{u}_t =
\mu(g(t),\tau)\tilde{u}_t,$$ since $f_t$ is a minimizer for
$\mathcal{W}$. Since $\mu(g(t),\tau)$ is uniformly bounded, as in the
previous section we can get that $C^{2,\alpha}$ norms of $\tilde{u}_t$
are uniformly bounded.  This implies that $C^{2,\alpha}$ norms of
$u_t$ are uniformly bounded.  Before we proceed with further
discussion notice the following.
\begin{remark}
$\int_M (4\pi\tau)^{-\frac{n}{2}}e^{-f_t(s)} dV_{g(s)} = 1$. This is a
simple consequence of the fact that $\int_M
(4\pi\tau)^{-\frac{n}{2}}e^{-f_t} dV_{g(t)} = 1$, since $f_t$ is a
minimizer for $\mathcal{W}$ with respect to $g(t)$, and the following
backward parabolic equation
$$\frac{d}{ds}f_t(s) = -\Delta f_t(s) + |\nabla f_t(s)|^2 - R +
\frac{n}{2\tau}.$$
Namely,
\begin{eqnarray*}
\frac{d}{ds}(\int_M e^{-f_t(s)} dV_{g(s)}) &=& \int_M
e^{-f_t(s)}(\Delta f_t(s) - |\nabla f_t(s)|^2 + R - \frac{n}{2\tau}
- R + \frac{n}{2\tau}) dV_{g(s)}\\
&=& \int_M \Delta (e^{-f_t(s)} dV_{g(s)}) = 0
\end{eqnarray*}
\end{remark}

Since $\log$ is a concave function and
$\tilde{u}_t(s)^2(4\pi\tau)^{-n/2}dV$ is a probability measure, we
have by Jensen and Sobolev inequalities
\begin{eqnarray*}
\int_M\tilde{u}_t(s)^2\log \tilde{u}_t(s)^2(4\pi\tau)^{-n/2}dV &=& 
\frac{n-2}{2}\int_M\tilde{u}_t(s)^2
\log \tilde{u}_t(s)^{4/(n-2)}(4\pi\tau)^{-n/2}dV \\
&\leq& \frac{n-2}{2}\log\int_M\tilde{u}_t(s)^{2n/(n-2)}(4\pi\tau)^{-n/2}dV \\
&\leq& \frac{n-2}{2}\log [C\int_M(|\nabla \tilde{u}_t(s)|^2+
\tilde{u}_t(s)^2)dV]^
{(n-2)/n} + \\
&+& \frac{n-2}{2}\log(4\pi\tau)^{-n/2} \\
&=& \frac{n}{2}\log C\int_M\tau(|\nabla \tilde{u}_t(s)|^2+\tilde{u}_t(s)^2)
(4\pi\tau)^{-n/2}dV.
\end{eqnarray*}
This inequality shows that 
\begin{equation}
\tau\int_M|\nabla \tilde{u}_t(s)|^2(4\pi\tau)^{-n/2}dV\leq C.
\end{equation}
The constant $C$ does not depend either on $t$ or $s\in[0,t]$. To 
conclude, we have the following estimates
$$\int_M |\tilde{u}_t (s)|^2 (4\pi \tau)^{-\frac{n}{2}}dV_s \le C_1$$
$$\tau (4\pi \tau)^{-\frac{n}{2}}\int_M |\nabla _s \tilde{u}_t(s)|^2 dV_s
\le C_2,$$
that is we have that $|\tilde{u}_t|_{W_{1,2}} \le C$ for a uniform
constant $C$.

Take a sequence $t_i\to\infty$. There exists a subsequence such that
$g(t_i+t)\to h(t)$ when $i\to \infty$, where $h(t)$ is a Ricci flow on
$M$. This follows from Hamilton's compactness theorem
(\cite{hamilton1995}). Fix $A>0$. $f_t$ will be a minimizer for
$\mathcal{W}$ with respect to $g(t)$, which we flow backward, for
every $t$. Let $s\in [0,A]$.

\begin{lemma}
\label{lemma-lemma_lower_bound}
For every $A > 0$ there exists $\delta = \delta(A) > 0$ such that
$u_{t+A}(t+s) \ge \delta > 0$ for all $t$ and all $s\in [0,A]$.
\end{lemma}

\begin{proof}
Assume that the statement of the lemma is not true. In that case there
would exist a sequence $s_i$ such that $\min_M u_{s_i+A}(s_i + a_i)
\to 0$ as $i\to\infty$, for some $a_i\in[0,A]$. Consider the equation
\begin{eqnarray*}
\frac{d}{dt}u_{s_i+A}(s_i+t) &=& -\Delta u_{s_i+A}(s_i+t) + 
(R-\frac{n}{2\tau})u_{s_i+A}(s_i+t),\\
u_{s_i+A}(s_i+A) &=& u_{s_i+A},
\end{eqnarray*}
for $t\in[0,A]$.  Let $\hat{u}_i(s_i+t) = \min_M
u_{s_i+A}(s_i+t)$. Then $\Delta \hat{u}_{s_i+A}(s_i+t) \ge 0$ and
$$\frac{d}{dt}\hat{u}_i(s_i+t) \le C\hat{u}_i(s_i+t),$$ 
where $C$ is a uniform constant. 
If we integrate it with respect to $t$, we get
$$\hat{u}_i(s_i+A) \le e^{CA} \hat{u}_i(s_i+t).$$ 
Since $\hat{u}_i(s_i+A) = \min_M u_{s_i+A}$ and since by Lemma
\ref{lemma-lower_bound} we know that there exists a constant $\delta$
such that $u_{s_i+A} \ge \delta > 0$, we have that $u_{s_i+A}(s_i+t)
\ge \delta(A) > 0$ for all $i$ and all $t\in [0,A]$.  This contradicts
our assumption that $\hat{u}_i(s_i + a_i) \to 0$ as $i\to\infty$.
\end{proof}

\begin{lemma}
For every $A > 0$ there exists $C(A)$ such that
\begin{enumerate}
\item
$\int_M u_{t}(s)^2 dV_{g(s)} \le C(A).$
\item
$\int_M |\nabla u_{t}(s)|^2 dV_{g(s)} \le C(A),$
\end{enumerate}
for all $t \ge A$, $s\in[t-A,t]$.
\end{lemma}

\begin{proof}
We will consider the equation
\begin{eqnarray*}
\frac{d}{ds}u_t(s) &=& -\Delta u_t(s) + (R-\frac{n}{2\tau})u_t(s) \\
u_t(t) &=& u_t,
\end{eqnarray*}
where $u_t = e^{-f_t}$ and $f_t$ is a minimizer for $\mathcal{W}$ with
respect to metric $g(t)$. Let $\hat{u}_t(s) = \max_M u_t(s)$. Then

$$\frac{d}{ds}\hat{u}_t(s) \ge -C\hat{u}_t(s),$$
where $C > 0$ is a uniform constant that does not depend either
on $s$ or $t$, but on the uniform bounds on geometries $g(t)$.
If we integrate it with respect to $s$ we get

$$\hat{u}_t = \hat{u}_t(t) \ge e^{-CA}\hat{u}_t(s),$$ for any $s\in
[t-A,t]$. On the other hand, we have already proved in the previous
section that $C^{2,\alpha}$ norms of $u_t$ are uniformly bounded in
$t\in [0,\infty)$.  Therefore we get that $0 \le u_t(s) \le C(A)$ on
$M$ for all $t\in [A,\infty)$ and all $s\in [t-A,t]$. Now we
immediately get part $1$ of our claim.  For part $2$ notice that

$$\int_M |\nabla u_t(s)|^2 dV_{g(s)} = 4\int_M u_t(s)|\nabla
\tilde{u}_t(s)|^2 dV_{g(s)} \le \tilde{C}(A),$$
since $\int_M |\nabla\tilde{u}_t(s)|^2$ is uniformly bounded for all 
$t \ge A$ and $s\in [t-A,t]$. 
\end{proof}

The previous two lemmas tell us that in order to find the uniform
estimates on $f_{t_i+A}(t_i+s)$ for $s\in [0,A]$, it is enough to find
the uniform $C^{k,\alpha}$ estimates on $u_{t_i+A}(t_i+s)$. Our main
goal in this section is to prove the following theorem.

\begin{theorem}
\label{theorem-uniform_estimate}
Under the assumptions of the main theorem, with the notations as
above, for every $A>0$ there exists a uniform constant $C$, depending
on $A$ such that $|u_{t}(s)|_{C^{2,\alpha}}\le C$ for all $t\ge A$,
$\forall s\in[t-A,t]$.
\end{theorem}

\begin{proof}
Consider the equation
$$\frac{d}{ds}u_t(s) = -\Delta u_t(s) + (R(s) - \frac{n}{2\tau})u_t(s),$$
for $t\in [A,\infty)$ and $s\in[t-A,t]$. All our further estimates will depend 
on $A$. We will use $C$ to denote different absolute constants that
depend on $A$ and the uniform bounds on our geometries $g(t)$.
Denote by $h = h_t(s) = (-\frac{n}{2\tau} + R(s))u_t(s)$. Omit the
subscript $t$.
$$\frac{d}{ds}u + \Delta u = h.$$
\begin{equation}
\label{equation-integral}
\int_M h^2 = \int_M (\frac{d}{ds}u)^2 + 2\int_M \frac{d}{ds}u\Delta u +
\int_M (\Delta u)^2,
\end{equation}
where we should keep in mind that the metric depends on $s$.
\begin{eqnarray}
\label{equation-equation_integrals}
\int_M \frac{d}{ds}u \Delta u &=& -\int_M g^{ij}\nabla_i
(\frac{d}{ds}u)\nabla_j u dV_s \\ \nonumber 
&=& -\frac{1}{2}\frac{d}{ds}\int_M|\nabla u|^2 dV_s - \int_M|\nabla
u|^2(\frac{n}{2\tau} - R)dV_s + \nonumber \\
&+& \int_M g^{pi}g^{qj}D_iu D_ju (2R_{pq} - 
\frac{1}{2\tau}g_{pq})dV_s,\nonumber
\end{eqnarray}
where the second term on the right hand side of
(\ref{equation-equation_integrals}) comes from taking the derivative
of the volume element and the third term appears from taking the
derivative of $g^{ij}$. Denote the former one by $J_1$ and the latter
one by $J_2$.

\begin{eqnarray*}
\int_M (\Delta u)^2 & = & \int_M g^{ij} D_iD_ju g^{kl}D_kD_lu \\
& = & - \int_M g^{ij}g^{kl}D_juD_iD_kD_lu \\
&=& - \int_Mg^{ij}g^{kl}D_juD_kD_iD_lu + \int_M
g^{ij}g^{kl}D_juR^l_{iks}D_su \\
&=& I + \int_Mg^{ij}g^{kl}D_kD_juD_iD_lu \\
&=& I + \int_M |\nabla^2u|^2,
\end{eqnarray*}
where $I = \int_M g^{ij}g^{kl}D_juR^l_{iks}D_su$.
Let $l\in (t-A,t)$ where $A>0$. Integrating the equation
(\ref{equation-integral}) in $s$, from $l$ to $t$ gives		
\begin{eqnarray*}							
\int_{l}^t(\int_M (\frac{d}{ds}u)^2 dV_s)ds + \int_M|\nabla u|^2 dV_s
|_{s=l} + \int_{l}^t\int_M|\nabla^2 u|^2dV_s ds \\ = \int_l^t\int_M h^2
+ \int_M|\nabla u|^2dV_s |_{s=t} + \int_l^t (2J_1 + 2J_2 + I).
\end{eqnarray*}

$$\int_l^t J_1 \le AC\sup_{s\in(t-A,t)}\int_M|\nabla u|^2 dV_s \le
\tilde{C},$$
for every $t$. 
Similarly we get estimates for $J_2$ and $I$. From all these estimates
we can conclude the following
\begin{equation}
\label{equation-estimate1}
\int_{t-A}^t\int_M (\frac{d}{ds}u_t(s))^2 dV_s ds \le C.
\end{equation}

\begin{equation}
\label{equation-estimate2}
\int_{t-A}^t\int_M|\nabla^2u_t(s)|^2dV_s ds \le C.
\end{equation}
\be
\label{equation-estimate1_2}
\sup_{s\in(t-A,t)}\int_M|\nabla u|^2 dV_s \le C, \ee 
where $C = C(A)$.
Let $\tilde{u}_t=\frac{d}{ds}u_t(s)$ (we will not confuse this
$\tilde{u}_t$ with one defined at the beginning of this section). Omit
the subscript $t$.
$$\frac{d}{ds}\tilde{u} = -D_s\Delta_su + \frac{d}{ds}[(R -
\frac{n}{2\tau})u].$$ 
Multiply the equation by $\tilde{u}$ and integrate it along $M$.
\begin{eqnarray*}
\frac{1}{2}\frac{d}{ds}\int_M|\frac{d}{ds}u|^2dV_s &=&
-\int_M\frac{d}{ds}(g(s)^{ij}D_iD_ju)\tilde{u} +
\int_M(\frac{d}{ds}(R - \frac{n}{2\tau}))u\tilde{u} + \frac{1}{2} 
\int_M (R - \frac{n}{2\tau})|\frac{d}{ds}u|^2 dV_s\\ 
&=& 2\int_M (-R_{pq} + \frac{1}{2\tau}g_{pq})g^{pi}(s)g^{qj}(s)D_iD_ju
\tilde{u} - \int_M g(s)^{ij}D_iD_j(\frac{d}{ds}u)\tilde{u} + \\
&+& \int_M(\frac{d}{ds}(R - \frac{n}{2\tau}))u\tilde{u} +  
\int_M g^{jk}(\frac{d}{dt}\Gamma_{ij}^k)\frac{\partial u}
{\partial x_k}\tilde{u} 
+ \frac{1}{2}\int_M (R - \frac{n}{2\tau})|\frac{d}{ds}u|^2dV_s.
\end{eqnarray*}
Since $\int_M g(s)^{ij}D_iD_j\frac{d}{ds}u\tilde{u} = -\int_M |\nabla_s
(\frac{d}{ds}u)|^2$ and since we are on the Ricci flow, metrics $g(s)$
are uniformly bounded, after applying Cauchy-Schwartz inequality and
using the uniform boundedness of the curvature operator, we get
\begin{eqnarray*}
& &\int_{t-A}^t\int_M |\nabla (\frac{d}{ds}u)|^2dV_s ds +
\sup_{s\in(t-A,t)}\int_M |\frac{d}{ds}u|^2 \le\\
&\le& C\int_{t-A}^t\int_M|\frac{d}{ds}u|^2dV_s ds +
C\int_{t-A}^t\int_M|\nabla^2u|^2dV_s ds + \\
&+& \int_M |\frac{d}{ds}u|^2 dV_s |_{s=t} + C\int_M |\nabla u|^2.
\end{eqnarray*}
$\int_M |\frac{d}{ds}u(s)|^2dV_s |_{s=t} \le C(\int_M |\Delta u_t|^2 +
\int_M h(t)^2)$ where $h(s)=(\frac{n}{2\tau} - R(s))u(s)$. Since $u_t
= e^{-f_t}$, where $f_t$ are the minimizers for $\mathcal{W}$, like in
the previous section we can conclude that $u_t \in W^{k,p}$, with
uniform bounds on $W^{k,p}$ norms (these bounds depend on $k$) and
therefore, $\int_M|\frac{d}{ds}u(s)dV_s |_{s=t}$ are uniformly bounded
in $t$.  This estimate together with estimates
(\ref{equation-estimate1}) and (\ref{equation-estimate2}) gives that
\begin{equation}
\label{equation-estimate3}
\int_{t-A}^t\int_M|\nabla (\frac{d}{ds}u)|^2dV_s ds \le C.
\end{equation}
\begin{equation}
\label{equation-estimate4}
\sup_{s\in (t-A,t)}\int_M |\frac{d}{ds}u|^2 \le C.
\end{equation}
If $\tilde{u} = \frac{d}{ds}u$ and $\tilde{h}=\frac{d}{ds}h$ then:
$$\frac{d}{ds}\tilde{u} = -D_s\Delta u + \tilde{h}.$$
\begin{eqnarray*}
D_s\Delta u = \frac{d}{ds}(g(s)^{ij}D_iD_j u) &=&
g(s)^{ip}g(s)^{jq}(\frac{1}{\tau}g_{pq}-2R_{pq})D_iD_j u +
g(s)^{ij}D_iD_j\tilde{u} \\ 
&+& g(s)^{ij}\frac{d}{ds}(\Gamma_{ij}^k)D_k u.
\end{eqnarray*}
\begin{eqnarray}
\label{equation-equation_for_tilde}
H &=& \tilde{h} - g^{ip}g^{jq}(\frac{1}{\tau}g_{pq}-2R_{pq})D_iD_j u -
g(s)^{ij}\frac{d}{ds}(\Gamma_{ij}^k)D_k u \\
&=& \frac{d}{ds}\tilde{u} + \Delta \tilde{u}. \nonumber
\end{eqnarray}
All the estimates that we have got so far tell that
$\int_{t-A}^t\int_M H^2$ is uniformly bounded in $t$. The analogous
estimates to the estimates (\ref{equation-estimate1}),
(\ref{equation-estimate2}) and (\ref{equation-estimate1_2}) for $u$, we
can get for $\frac{d}{ds}u$ (by using the evolution equation for
$\frac{d}{ds}u$ and all the estimates that we have got so far by
analyzing the evolution equation for $u$).
\begin{equation}
\label{equation-estimate5}
\int_{t-A}^t\int_M(|\nabla^2(\frac{d}{ds}u)|^2 dV_s ds \le C.
\end{equation}

\begin{equation}
\label{equation-estimate6}
\int_{t-A}^t\int_M(\frac{d^2}{ds^2}u)^2dV_sds \le C.
\end{equation}

\begin{equation}
\label{equation-estimate7}
\sup_{s\in(t-A,t)}\int_M |\nabla (\frac{d}{ds}u)|^2 dV_s \le C.
\end{equation}
To obtain these estimates we have used the fact that
$$\int_M |\nabla\frac{d}{ds}u|^2dV_{g(s)}|_{s=t} \le
C(\int_M|\nabla\Delta u_t|^2 + \int_M|\nabla(R-\frac{n}{2\tau})u_t|^2,$$
where the right hand side is uniformly bounded in $t$, since
$u_t = e^{-f_t}$ and $f_t$ are the minimizers for $\mathcal{W}$.

By standard regularity theory, considering $\Delta u_t(s) =
-\frac{d}{ds}u_t(s) + h_t(s)$ as an elliptic equation whose right hand
side has uniformly bounded $W^{1,2}$ norms for $s\in (t-A,t)$ and all
$t\ge A$, we have that $|u_t(s)|_{W^{3,2}} \le C$, for a uniform
constant $C$ that depends on $A$. Take a derivative in $s$ of the
equation $\frac{d}{ds}\tilde{u} = -\Delta\tilde{u} + H$, with
$\tilde{u} =
\frac{d}{ds}u$. Denote by $\bar{u} =
\frac{d}{ds}\tilde{u}$. By using the estimates that we have got for
$\tilde{u}$ it is easy to conclude that $\bar{u}$ satisfies the
equation
$$\frac{d}{ds}\bar{u} = -\Delta\bar{u} + H_1,
$$ where $H_1 = \frac{d}{ds}H + g^{ip}g^{jq}(-2R_{pq}+\frac{1}{\tau})
D_iD_j \tilde{u} + g(s)^{ij}\frac{d}{ds}(\Gamma_{ij}^k)D_k \tilde{u} $ and
$\int_{t-A}^t\int_M H_1^2 dV_{g(s)}ds$ is uniformly bounded in $t$.
As in the case of the previous estimates we can conclude that

$$\sup_{s\in(t-A,t)}\int_M|\frac{d}{ds}\tilde{u}|^2dV_s \le C,$$
$$\sup_{s\in(t-A,t)}\int_M|\nabla(\frac{d}{ds}\tilde{u})|^2 dV_s \le
C.$$ 
By regularity theory applied to the equation $\Delta\tilde{u} =
-\frac{d}{ds}\tilde{u} + H$, we can get that
$\frac{d}{ds}u_t(s)$ has uniformly bounded $W^{3,2}$ norms. If we
go back to the parabolic equation for $u_t(s)$ we can get that
$|u_t(s)|_{W^{5,2}} \le C$ for all $t\ge A$  and all $s\in (t-A,t)$. 
Continuing this process by taking more and more derivatives in $t$ 
of our original parabolic equation we can
conclude that $W^{p,2}$ norms of $u_t(s)$ are uniformly
bounded for every $p$, by the constants that depend on $A$ and
$p$. Sobolev embedding theorem now gives that
all $C^{k,\alpha}$ norms of $u_t(s)$ are uniformly bounded
for all $t > A$ and all $s\in [t-A,t]$, by  constants that depend
on $A$ and $k$.
\end{proof}

Combining Theorem \ref{theorem-uniform_estimate} and Lemma
\ref{lemma-lemma_lower_bound}, we get that for every $A$ there exist
constants $C_k = C(k,A)$ such that $|f_t(s)|_{C^{k,\alpha}} \le C_k$,
for all $t \ge A$ and all $s\in [t-A,t]$.

\end{subsection}

\begin{subsection}{Ricci soliton in the limit}
\label{section:soliton}

In this subsection we want to finish the proof of Theorem
\ref{theorem-theorem_main}. 

We have uniform curvature and diameter bounds for our flow $g(t)$. We
have already proved that we also have a volume noncollapsing condition
along the flow, for all times $t \ge 0$. This gives a uniform lower
bound on the injectivity radii. Hamilton's compactness theorem
(modified to the case of our flow) gives that for every sequence
$t_i\to\infty$ there exists a subsequence so that $g(t_i+t)\to h(t)$
uniformly on compact subsets of $M\times [0,\infty)$ and that $h(t)$
is a solution to the Ricci flow (\ref{equation-equation_our_flow}). We
will show below that for each $t$, $h(t)$ satisfies actually a Ricci
soliton equation with the Hessian of function $f_h(t)$ involved, where
$f_h(t)$ is a smooth one parameter family of functions. We will now
see how we get the functions $f_h(t)$, using the estimates on $f_t(s)$
from the previous subsection and Perelman's monotonicity formula.

Take any $t$ and let $f_t$ be a function so that $\mu(g(t),\tau) =
\mathcal{W}(g(t),f_t,\tau)$. Flow $f_t$ backward. Fix $A > 0$. Then:
$$I(t) = \mathcal{W}(g(t+A),f_{t+A},\tau) -
\mathcal{W}(g(t),f_{t+A}(t),\tau) \le \mu(g(t+A),\tau)-\mu(g(t),\tau)
\to 0 (t\to \infty).$$
$$0 \le I(t) = \int_0^A \frac{d}{du}W(g(t+s),f_{t+A}(t+s),\tau)ds \to
0,$$ as $t\to\infty$.  
We will consider $u_{t_i+A}(t_i+s)$ where $s\in
[0,A]$.  We will divide the proof of the theorem in a few steps.

\begin{step}
\label{step-step_soliton}
$\forall A>0$, $\lim _{i\to \infty}
\frac{d}{du}W(g(s+t_i),f_{t_i + A}(s+t_i),\tau) = 0$ for almost all
$s\in [0,A]$.
\end{step}

\begin{proof}
$I(t_i)\to 0$ by Claim \ref{claim-step0}. On the other hand
$$I(t_i) = \mathcal{W}(g(t_i+A),f_{t_i+A},\tau) -
\mathcal{W}(g(t_i),f_{t_i+A}(t_i),\tau) =
\int_0^A\frac{d}{du}W(g(t_i+s),f_{t_i+A}(t_i+s),\tau)ds.$$ Since by
Perelman's monotonicity formula
$\frac{d}{du}W(g(t_i+s),f_{t_i+A}(t_i+s),\tau) \ge 0$, we have that
$\lim_{i\to\infty}\frac{d}{du}W(g(t_i+s),f_{t_i+}(t_i+s),\tau) = 0$
for almost all $s\in [0,A]$, for

$$\int_0^A\lim_{i\to\infty}\frac{d}{du}W(g(t_i+s),f_{t_i+A}(t_i+s),\tau)ds
\le \lim_{i\to\infty} I(t_i),$$
by Fatuous lemma.
\end{proof}

\begin{step}
\label{step-step_bound_on_s_derivatives}
$|\tilde{u}_t(s)|_{C^{2,\alpha}} \le C$, $\forall t$, where
$\tilde{u}_t(s)=\frac{d}{ds}u_t(s)$.
\end{step}

\begin{proof}
Following the notation of the previous subsection, we get that:
$$\frac{d}{ds}\tilde{u}_t(s) = -\Delta \tilde{u}_t(s) + H_t(s),$$
where $H_t(s) = \frac{d}{ds}h_t(s) + g^{ip}g^{jq}(\frac{1}{\tau}g_{pq}
-2R_{pq})D_iD_ju + g^{ij}\frac{d}{ds}(\Gamma_{ij}^k)D_k u$.
$$\tilde{u}_t(t) = \frac{d}{ds}u_t(s) = -\Delta u_t +
(-\frac{n}{2\tau} +R)u_t.$$ 
In the previous subsection we have proved
that there exist a uniform lower and an upper bound on $u_t(s)$ and
that $|u_t(s)|_{W^{3,p}}\le C(p,A)$ for all $t\ge A$ and all
$s\in[t-A,t]$. Similarly we can get that $|u_t(s)|_{W^{k,p}}\le
C(k,p,A)$ and therefore $|\tilde{u}_t(s)|_{W^{k-2,p}}\le C(k,p,A)$,
$\forall t\ge A$ and all $s\in[t-A,t]$. We can get that
$|\tilde{u}_t(s)|_{C^{2,\alpha}}\le C$, for all $t\ge A$ and $\forall
s\in [t-A,t]$.  We can extend this to all higher order time
derivatives of $u_t(s)$.
\end{proof}

\begin{step}
\label{step-step_partial}
For every $A > 0$ there exists a subsequence $t_i$, so that the limit
metric $h(s)$ of a sequence $g(t_i+s)$ is a Ricci soliton for $s\in
[0,A]$.
\end{step}

\begin{proof}
By step \ref{step-step_soliton} we have that
$$\lim_{i\to \infty}R_{jk}(t_i+s) + \nabla_j\nabla_k f_{t_i+A}(t_i +
s) -\frac{1}{2\tau}g_{jk}(t_i+s) = 0,$$ for almost all $s\in[0,A]$ and
almost all $x\in M$, since
$$\frac{d}{ds}\mathcal{W}(g(t_i+s),f_{t_i+A}(t_i+s),\tau) =
(4\pi\tau)^{-\frac{n}{2}}\int_M2\tau|R_{jk} +
\nabla_jf_{t_i+A}\nabla_kf_{t_i+A} - \frac{1}{2\tau}g_{jk}|^2
dV_{g(t_i+s)}.$$ By Lemma \ref{lemma-lemma_lower_bound} and Theorem
\ref{theorem-uniform_estimate}, we have that $0 < C_1 \le
|u_{t_i+A}(s+t_i)| \le C_2$ for all $i \ge i_0$ and all $s\in[0,A]$,
for some constants $C_1$ and $C_2$ that depend on $A$. By step
\ref{step-step_bound_on_s_derivatives} and Theorem
\ref{theorem-uniform_estimate} we can find a subsequence, say
$\{t_i\}$ such that $f_{t_i+A}(t_i+s)$ converges in $C^{2,\alpha}$
norm to $\tilde{f}_A(s)$ for all $s\in [0,A]$ and all $x\in M$. More
precisely, for a countable dense subset $\{s_j\}$ of $[0,A]$ there
exists a subsequence so that $f_{t_i+A}(t_i+s_j)$ converges in
$C^{2,\alpha}$ norm to $\tilde{f}_A(s_j)$ on $M$. For any $s\in [0,A]$
there exists a subsequence $t_{i_k}$ so that
$f_{t_{i_k}+A}(t_{i_k}+s)$ converges to $\tilde{f}_A(s)$ in
$C^{2,\alpha}$ norm. We want to show that actually
$f_{t_i+A}(t_i+s)\stackrel{C^{2,\alpha}}{\to}\tilde{f}_A(s)$. For that
we use the fact that $\frac{d}{ds}f_{t_i+A}(t_i+s)$ is uniformly
bounded in $C^{2,\alpha}$ norm, and therefore

$$|\tilde{f}_A(s) - \tilde{f}_A(s_0)|_{C^{2,\alpha}} < \epsilon,$$
for some small $\epsilon > 0$ and some $s_0\in \{s_j\}$ that is 
sufficiently close to $s$. We also have
$$|\tilde{f}_A(s_0) - f_{t_i+A}(t_i+s_0)|_{C^{2,\alpha}} < \epsilon,$$
for $i\ge i_0$ and
$$|f_{t_i+A}(t_i+s) - f_{t_i+A}(t_i+s_0)|_{C^{2,\alpha}} < \epsilon,$$
since $|\frac{d}{ds}f_{t_i+A}(t_i+s)|_{C^{2,\alpha}} \le C(A)$, for
all $i\ge i_0$ and all $s\in [0,A]$. By triangle inequality, we now
get that for every $\epsilon > 0$ there exists $i_0$ so that

$$|\tilde{f}_A(s) - f_{t_i+A}(t_i+s)|_{C^{2,\alpha}} < 3\epsilon,$$
for all $i\ge i_0$ and all $s\in [0,A]$.

$f_{t_i+A}(t_i+s)$ converges in
$C^{2,\alpha}$ norm on $M$ to $\tilde{f}_A(s)$, for all
$s\in[0,A]$. Finally, we get that
\be
\label{equation-equation_ricci_soliton}
R_{jk} + \nabla_j\nabla_k \tilde{f}_A(s) - \frac{1}{2\tau}h_{jk}(s) =
0, \ee 
for all $s\in [0,A]$, and for almost all $x\in M$. Because of
the continuity it will hold for all $x\in M$. Since $h(s)$ is a Ricci
flow, all covariant derivatives of $h$ and the covariant derivatives
of a curvature operator are uniformly bounded, and therefore
$|\nabla^p\tilde{f}_A(s)|\le C(p), \:\:\: \forall s\in [0,A]$ and all
$p\ge 2$. Also we have that $|\frac{d^p}{ds^k}\nabla ^p
\tilde{f}_A(s)| \le C(p,k)$ where $C(p,k)$ does not depend on $A$, for
$p \ge 2$.
\end{proof}

\begin{step}
We can glue all the functions $\tilde{f}_A$ that we get for different
values of $A$, to get a function $f_h(s)$ defined on
$M\times[0,\infty)$, which defines our metric $h(s)$ as a soliton type
solution for all times $s\ge 0$.
\end{step}

\begin{proof}
Take any increasing sequence $A_j\to\infty$. For every $A_j$, by the
previous step we can extract a subsequence $t_i$ so that
$f_{t_i+A_j}(t_i+s)\stackrel{C^{2,\alpha}}{\to}\tilde{f}_{A_j}(s)$ for
all $s\in[0,A_j]$. Diagonalization procedure gives a subsequence so
that $f_{t_i+A_j}(s)\stackrel{C^{2,\alpha}}{\to}\tilde{f}_{A_j}(s)$
for all $j$ and all $s\in [0,A_j]$. For this subsequence $t_i$ we have
that $g(t_i+t)\to h(t)$, uniformly on compact subsets of
$\times[0,\infty)$.  Compare the functions $\tilde{f}_{A_j}$ and
$\tilde{f}_{A_k}$ for $j < k$, on the interval $[0,A_j]$. We know that
they both satisfy

$$\Delta_{h(s)}\tilde{f}_{A_r} + R(h(s)) - \frac{n}{2\tau} = 0,$$ and
therefore $\Delta_{h(s)}(\tilde{f}_{A_j} - \tilde{f}_{A_k}) = 0$.
Since $M$ is compact, this implies that $\tilde{f}_{A_k}(s) =
\tilde{f}_{A_j}(s) + c_{A_j}^{A_k}(s)$, for $s\in[0,A_j]$, where $
c_{A_j}^{A_k}(s)$ is a constant function for every $s\in [0,A_j]$.  On
the other hand, because of the integral normalization condition, we
have

$$(4\pi\tau)^{-\frac{n}{2}}\int_M e^{-\tilde{f}_{A_j}(s)}dV_{h(s)} = 1,$$
$$(4\pi\tau)^{-\frac{n}{2}}\int_M e^{-\tilde{f}_{A_k}(s)}dV_{h(s)} = 1
= e^{-c_{A_j}^{A_k}(s)}(4\pi\tau)^{-\frac{n}{2}} \int_M
e^{-\tilde{f}_{A_j}(s)}dV_{h(s)},$$ which implies that
$c_{A_j}^{A_k}(s) = 0$ for all $s\in [0,A_j]$ and all $k \ge
j$. Therefore $\tilde{f}_{A_j}(s) = \tilde{f}_{A_k}(s)$ for all $s\in
[0,A_j]$.  Define a function $f_h(s)$ in the following way. Let
$f_h(s) = \tilde{f}_{A_j}(s)$, for all $s\in[0,A_j]$ and all
$A_j\to\infty$. $f_h(s)$ is a well defined function because of the
previous discussion. We also have that

\begin{equation}
\label{equation-equation_soliton_help}
R(h(s))_{pq} + \nabla_p\nabla_q f_h(s) - \frac{1}{2\tau}h(s)_{pq} = 0,
\end{equation}
holds for all $s\in [0,\infty)$. The definition of $f_h(s)$ does not
depend on a choice of an increasing sequence $A_j$. Namely, if $B_j$
were another increasing sequence and if $f_{h'}(s)$ were functions
defined using the sequences $B_j$ and $t_i$ ($t_i$ is the same
sequence as above), then at each time both functions $f_h(s)$ and
$f_{h'}(s)$ would satisfy the same equation
(\ref{equation-equation_soliton_help}) and the same integral
normalization condition. Therefore $f_h(s) = f_{h'}(s)$ for all $s\in
[0,\infty)$.
\end{proof}

\end{subsection}

\begin{subsection}{Some properties of the limit solitons}

Let $t_i$ be any sequence converging to infinity. Then as we have seen
earlier, there exists a subsequence such that $g(t_i+s)\to h(s)$,
where $h(s)$ is a Ricci soliton. Let $\hat{R}(h(t))=\min R(h(t))$. We
will first state a theorem that R. Hamilton proved in his paper
\cite{hamilton1999}.

\begin{theorem}[Hamiton]
\label{theorem-theorem_hamilton}
Under the normalized Ricci flow, whenever $\hat{R}\le 0$, it is
increasing, whereas if ever $\hat{R}\ge 0$ it remains so forever.
\end{theorem}

We will use the proof of Theorem \ref{theorem-theorem_hamilton} to
prove the following lemma.

\begin{lemma}
Under the assumptions of Theorem \ref{theorem-theorem_main},
$\hat{R}(h(t))\ge 0$, $\forall t$, for the limit metric $h(t)$ of any
sequence of metrics $g(t_i)$, where $g(t)$ is a solution of

\[\frac{d}{dt}g_{jk} = -2R_{jk}(g(t)) + \frac{1}{\tau}g_{jk}(t).\]
\end{lemma}

\begin{proof}
Assume that there exists $t_0$ such that $\hat{R}(h(t_0))<0$. Without loss
of generality assume that $t_0=0$. Since $g(t_i)\to h(0)$ as $i\to
\infty$, there exists $i_0$, so that for all $i\ge i_0$
$\hat{R}(g(t_i))<0$. The evolution equation for $R$ is
$$\frac{d}{dt}R = \Delta R + 2|\ric|^2 +
\frac{2}{n}R(R-\frac{n}{2\tau}).$$
This implies
$$\frac{d}{dt}\hat{R} \ge
\frac{2}{n}\hat{R}(\hat{R}-\frac{n}{2\tau}).$$ If $\hat{R} \le 0$,
then $\hat{R}$ is increasing (since $\frac{d}{dt}\hat{R} \ge 0$). If
$\hat{R} \ge 0$ at some time it can not go negative at later times.
If there existed $t>t_{i_0}$ such that $\hat{R}(g(t))\ge 0$, then
$\hat{R}\ge 0$ would remain so forever, for all $s\ge t$ and therefore
we could not have $\hat{R}(g(t_i))<0$ for $t_i>t$. That contradicts
the fact that $\hat{R}(g(t_i))<0$ for all $i\ge i_0$. Therefore
$\forall t\ge t_{i_0}$ we have that $\hat{R}(g(t))<0$.
$$\frac{d\hat{R}}{dt}\ge \frac{2}{n}\hat{R}(\hat{R}-\frac{n}{2\tau})
\ge 0,$$
for all $t$ big enough. That implies $\hat{R}$ is increasing and
therefore there exists $\lim_{t\to \infty}\hat{R}(g(t)) = -C \le
0$. Moreover $\hat{R}(h(s))=-C$ for all $s$. Since $\lim_{i\to
\infty}\hat{R}(g(t_i)) = \hat{R}(h(0)) < 0$, $C>0$. We also have that

$$\frac{d\hat{R}(h(s))}{ds}\ge
-\frac{2}{n}\hat{R}(h(s))(\frac{n}{2\tau}-\hat{R}(h(s))) =
\frac{2}{n}C(\frac{n}{2\tau}+C) \ge 0.$$ The left hand side of the
above inequality is zero and therefore we get that $C =
-\frac{n}{2\tau}$ or $C = 0$. Since $C > 0$, we get a contradiction.
Therefore $R(h(t))\ge 0$ for all $t$, what we wanted to prove.
\end{proof}

\begin{remark}
\label{remark-remark_soliton}
Let $(M,g)$ be a compact manifold and $g(t)$ be a Ricci flow on $M$.
Since
$$\frac{d}{dt}\mathcal{W} = \int_M 2\tau|R_{ij} + \nabla_i\nabla_j f -
\frac{1}{2\tau}g_{ij}|^2(4\pi\tau)^{-\frac{n}{2}} e^{-f} dV,$$
$\mathcal{W}(g,f,\tau) = const$ along the flow, if $g$ is a Ricci soliton
satisfying the equation
$$R_{ij} + \nabla_i\nabla_j f - \frac{1}{2\tau}g_{ij} = 0.$$
\end{remark}
 
Let $t_i\to \infty$ and $s_i\to \infty$ be two sequences such that
$g(t_i + t)\to h(t)$ and $g(s_i + t)\to h'(t)$ where $h(t)$ and
$h'(t)$ are 2 Ricci solitons on $M$ that have been constructed
earlier. We have proved that
$$R_{jk}(h) + \nabla_j\nabla_k f_h(t) - \frac{1}{2\tau}h_{jk} = 0,$$
$$R_{jk}(h') + \nabla_j\nabla_k f_{h'}(t) - \frac{1}{2\tau}h'_{jk} = 0,$$
where

$$f_h(t) = \lim_{j\to\infty}\lim_{i\to\infty} f_{A_j + t_i} (t_i + t),$$
$$f_{h'}(t) = \lim_{j\to\infty}\lim_{i\to\infty} f_{B_j + s_i} (s_i + t),$$
for some increasing sequences $A_j\to\infty$ and $B_j\to\infty$.
By Remark \ref{remark-remark_soliton} we know that
$\mathcal{W}(h(t),f_h(t),\tau) = C_1$ and $W(h'(t),f_{h'}(t),\tau) = C_2$ are
constant along the flows $h(t)$ and $h'(t)$ respectively.

\begin{lemma}
\label{lemma-lemma_W_constant}
$C_1 = C_2$, i.e. $\mathcal{W}(h(t),f_h(t),\tau)$ is a same constant
for all solitons $h(t)$ that arise as limits of sequences of metrics
of our original flow $g(t)$ (\ref{equation-equation_our_flow}) on a
compact manifold $M$.
\end{lemma}

\begin{proof}
\begin{eqnarray}
\label{equation-equation_W_one}
& &\mathcal{W}(g(t_i + t),f_{t_i + A_j}(t_i + t),\tau) - \mathcal{W}(g(s_i),
f_{s_i + B_j}(s_i),\tau) \le \nonumber\\ 
&\le& \mathcal{W}(g(t_i + A_j),f_{t_i + A_j}(t_i + A_j),\tau) - 
\mathcal{W}(g(s_i),f_{s_i}(s_i),\tau) = \nonumber \\ 
&=& \mu(g(t_i + A_j),\tau) - \mu(g(s_i),\tau) \to 0,
\end{eqnarray}
where we have used the fact that $\mathcal{W}(g(t),f(t),\tau)$
increases in $t$ along the flow (\ref{equation-equation_our_flow}) and
the fact that $f_{s_i}(s_i) = f_{s_i}$ is a minimizer for
$\mathcal{W}(g(s_i),f,\tau)$ over all $f$ belonging to a set
$\{f\:\:|\:\:\int_M (4\pi\tau)^{-\frac{n}{2}}e^{-f}
dV_{g(s_i)}\}$. Similarly,
\begin{eqnarray}
\label{equation-equation_W_two}
& &\mathcal{W}(g(t_i + t),f_{t_i + A_j}(t_i + t),\tau) - 
\mathcal{W}(g(s_i), f_{s_i + B_j}(s_i),\tau) \ge \nonumber\\
&\ge& \mathcal{W}(g(t_i + t),f_{t_i + t}(t_i + t),\tau) - 
\mathcal{W}(g(s_i + B_j),f_{s_i + B_j}(s_i + B_j),\tau) = \nonumber\\
&=& \mu(g(t_i + t),\tau) - \mu(g(s_i + B_j),\tau) \to 0,
\end{eqnarray}
when $i\to\infty$. From equations (\ref{equation-equation_W_one}) and
(\ref{equation-equation_W_two}), letting $i\to\infty$ we get

$$\mathcal{W}(h(t),\tilde{f}_{A_j}(t),\tau) - 
\mathcal{W}(h'(0),\tilde{f}'_{B_j}(0),\tau) \le 0.$$
$$W(h(t),\tilde{f}_{A_j}(t),\tau) - W(h'(0),\tilde{f}'_{B_j}(0),\tau) \ge 0.$$
Let $j\to\infty$ to get 

$$C_1 = \mathcal{W}(h(t),f_h(t),\tau) =
\mathcal{W}(h'(0),f_{h'}(0),\tau) = C_2.$$
\end{proof}

\begin{lemma}
\label{lemma-lemma_minimizers_solitons}
For every Ricci soliton $h(t)$ that arises as a limit of some sequence
of metrics of our original flow $g(t)$, the corresponding function
$f_h(t)$, that we have constructed before, is a minimizer for
Perelman's functional $\mathcal{W}$ with respect to a metric $h(t)$.
\end{lemma}

\begin{proof}
We will first proof the following claim.

\begin{claim}
\label{claim-claim_minimizer_one_sequence}
There exists a sequence $t_i\to\infty$ such that $g(t_i + t)\to h(t)$
as $i\to\infty$, where $h(t)$ is a Ricci soliton satisfying $R_{jk}(h)
+ \nabla_j\nabla_k f_h - \frac{1}{2\tau}h_{jk} = 0$ and $f_h(t)$ is a
minimizer for $\mathcal{W}(h(t),f,\tau)$.
\end{claim}

\begin{proof}[Proof of the Claim]
Let $H(t) = (4\pi\tau)^{-n/2}\int_M 2\tau|R_ij + \nabla_i\nabla_j f_t
- \frac{1}{2\tau}g_{ij}|^2 dt$, where $f_t$ is a function such that
$\mu(g(t),\tau) = W(g(t),f_t,\tau)$.  If we flow $f_t$ backward by
the equation

$$\frac{d}{dt}f = -\Delta f + |\nabla f|^2 - R + \frac{n}{2\tau},$$
starting at time $t$, for every $t > 0$ we get solutions
$f_t(s)$. Look at $F_t(s) = \mathcal{W}(g(s),f_t(s),\tau)$. We know that

$$\frac{d}{ds}F_t(s) = (4\pi\tau)^{-\frac{n}{2}}
\int_M 2\tau|R_{jk} + \nabla_j\nabla_k f_t(s) 
- \frac{1}{2\tau}g(s)_{jk}|^2 dV_{g(s)}.$$
$F_t(s)$ is a continuous function in $s\in[0,t]$ and $\lim_{s\to
t}\frac{d}{ds}F_t(s) = H(t)$. Therefore there exists a left derivative
of $F_t(s)$ at point $t$ and $(F_t)_{-}'(t) = H(t)$ for every $t > 0$.
Moreover, $g(t)$ and all the derivatives of $f_t$ up to the second
order are Lipshitz functions in $t$ (this follows from the estimates in
the previous subsections) and therefore

$$\mu(t) := \mu(g(t),\tau) = \inf_{\{f\:|\: \int_M
(4\pi\tau)^{-\frac{n}{2}}e^{-f} = 1\}}\mathcal{W}(g(t),f,\tau)$$
is a Lipshitz function in $t$ as well, i.e.  $k(t) = F_t(t) =
\mathcal{W}(g(t),f_t,\tau)$ is a Lipshitz function in $t$. This tells 
that $k(t)$ is differentiable in $t$, almost everywhere. Our discussion
then implies that $k'(t) = H(t)$ in a sense of distributions.

\begin{eqnarray}
\label{equation-equation_finite_H}
\int_{\delta}^{\infty}H(t) dt &=& \lim_{K\to\infty}\int_{\delta}^K
k'(t) dt \nonumber \\ &=& \lim_{K\to\infty} W(g(K),f_K,\tau) -
W(g(\delta),f_{\delta},\tau) \nonumber\\ &=&
\lim_{K\to\infty}(\mu(g(K),\tau) - \mu(g(\delta),\tau) \le C,
\end{eqnarray}
where $\delta > 0$ and $C$ is some uniform constant.
We have that $\int_{\delta}^{\infty} H(t) \le C$. This implies that
there exists a sequence $t_i\to\infty$ such that $H(t_i)\to 0$ as
$i\to\infty$, i.e.

$$\lim_{i\to\infty}(R_{jk} + \nabla_j\nabla_k f_{t_i} -
\frac{1}{2\tau}g_{jk})(t_i) = 0.$$ 
By what we have proved before,
after extracting a subsequence we can assume that $g(t_i)\to h(0)$
smoothly and $f_{t_i}\to \tilde{f}$ in $C^{2,\alpha}$ norm, where by
Theorem \ref{theorem-theorem_continuity_of_minimizers} $\tilde{f}$ is
a minimizer for $\mathcal{W}$ with respect to metric $h(0)$. Therefore,
\begin{equation}
\label{equation-equation_one}
R_{jk}(h(0)) + \nabla_j\nabla_k\tilde{f} - \frac{1}{2\tau}h_{jk}(0) = 0.
\end{equation}
On the other hand $g(t_i + t)\to h(t)$ as $i\to\infty$ 
where $h(t)$ is a Ricci soliton and

\begin{equation}
\label{equation-equation_two}
R_{jk}(h(t)) + \nabla_j\nabla_k f_h(t) - \frac{1}{2\tau} h_{jk}(t) = 0,
\end{equation}
where $f_h(t) = \lim_{j\to\infty}\lim_{i\to\infty}f_{t_i+A_j}(t_i+t)$,
for some sequence $A_j\to\infty$.
From equations (\ref{equation-equation_one}) and
(\ref{equation-equation_two}) we have that $\Delta (f_h(0) - \tilde{f})
= 0$, i.e. $f_h(0) = \tilde{f} + C$ for some constant $C$. We know
that $\int_M (4\pi\tau)^{-\frac{n}{2}}e^{-\tilde{f}} dV_{h(0)} = 1$, since
$\tilde{f}$ is a minimizer. From the construction of $f_h(t)$ it
follows that $\int_M (4\pi\tau)^{-\frac{n}{2}}e^{-f_h(0)} dV_{h(0)} =
1$ and therefore $\tilde{f} = f_h(0)$.
Since there exists a finite limit, $\lim_{t\to\infty}\mu(g(t),\tau)$, we
have that $\mu(h(0),\tau) = \mu(h(t),\tau)$ for all $t$. This implies 
that
\begin{eqnarray*}
\mu(h(t),\tau) &=& \mu(h(0),\tau) = \mathcal{W}(h(0),\tilde{f},\tau)\\
&=& \mathcal{W}(h(0),f_h(0),\tau) = \mathcal{W}(h(t),f_h(t),\tau),
\end{eqnarray*}
where we have used the fact that $\mathcal{W}$ is constant 
along a soliton. This means that $f_h(t)$ is a minimizer for 
$\mathcal{W}$ with respect to a metric $h(t)$, for every $t \ge 0$.
\end{proof}

To continue the proof of Lemma \ref{lemma-lemma_minimizers_solitons}
take any sequence $s_i\to\infty$. By a sequential convergence of our
original flow $g(t)$ to Ricci solitons, after extracting a subsequence
we may assume that $g(s_i + t)\to h'(t)$ as $i\to\infty$ where $h'(t)$
is a Ricci soliton. Take a soliton $h(t)$ with the properties as in
Claim \ref{claim-claim_minimizer_one_sequence}. From the convergence
of $\mu(g(t),\tau)$ we know that $\mu(h'(t),\tau) = \mu(h(s),\tau)$
for all $t$ and all $s$.

\begin{equation}
\label{equation-equation_mu_W}
\mu(h'(t),\tau) = \mu(h(s),\tau) = \mathcal{W}(h(s),f_h(s),\tau).
\end{equation}
By Lemma \ref{lemma-lemma_W_constant} we have that
$\mathcal{W}(h(s),f_h(s),\tau) = \mathcal{W}(h'(t),f_{h'}(t),\tau)$
for all $s$ and $t$. Combining this with 
(\ref{equation-equation_mu_W}) gives  that $\mu(h'(t),\tau) =
\mathcal{W}(h'(t),f_{h'}(t),\tau)$, i.e. $f_{h'}(t)$ is a minimizer for $h'(t)$
for every $t$.
\end{proof}

One useful property of the sequential soliton limits of our flow
(\ref{equation-equation_our_flow}) is that all limit solitons are the
solutions of the normalized flow equation

$$\frac{d}{dt}h_{ij} = -2R_{ij} + \frac{2}{n}r(h(t))h_{ij},$$
where $r(h(t)) = \frac{1}{\vol_{h(t)}M}\int_M R(h(t)) dV_{h(t)}$. In
the case of any of our soliton limits, we have that $R(h(t)) + \Delta
f_h(t) - \frac{n}{2\tau} = 0$ and therefore $r = r(h(t)) =
\frac{n}{2\tau}$ for all $t \ge 0$.

\begin{remark}
Let $t_i\to\infty$ and $g(t_i+t)\to h(t)$, where $h(t)$ is an Einstein
metric with an Einstein constant $\frac{1}{2\tau}$. If $\vol_{h'}(M) =
\vol_h(M)$, for any other limit soliton $h'$, then $h'$ is an Einstein
metric with the same Einstein constant $\frac{1}{2\tau}$. 
\end{remark}

\begin{proof}
The fact that $h$ is Einstein metric implies that
$\nabla_i\nabla_j f_h = -2R_{ij} + \frac{1}{\tau}h_{ij} = 0$, that is
$\Delta f_h = 0$. Since $M$ is compact, $f_h = C$ such that
$(4\pi\tau)^{-n/2}e^{-C}\vol_h(M) = 1$. An easy computation shows that
$\mu(h,\tau) = \mathcal{W}(h,C,\tau) = C - \frac{n}{2}$, and therefore
$\mu(h',\tau) = \mu(h,\tau) = C - \frac{n}{2}$. Then,
$(4\pi\tau)^{-n/2}e^{-C}\vol_{h'}(M) = 1$, implies that $f=C$ is a
minimizer for $\mathcal{W}$ with respect to $h'$ as well. This yields
$$\tau(2\Delta f - |\nabla f|^2 + R(h')) + f - n = C - \frac{n}{2},$$
that is
$$R(h') = \frac{n}{2\tau}.$$
From 
$$\Delta f_{h'} = \frac{n}{2\tau} - R(h') = 0,$$
we get that $f_{h'} = C$ and therefore 
$$R_{ij}(h') + \nabla_i\nabla_j f_{h'} - \frac{1}{2\tau}h'_{ij} = 0,$$
yields $R_{ij}(h') = \frac{1}{2\tau}h'_{ij}$.
\end{proof}

In the discussion that follows we will use Moser's weak maximum
principle. We will state it below, for a reader's convenience.

\begin{lemma}[Moser's weak maximum principle]
\label{lemma-lemma_Moser}
Let $g = g(t)$, $0 \le t < T$, be a smooth family of metrics, $b$ a
nonnegative constant and $f$ a nonnegative function on $M\times [0,
T)$ which satisfies the partial differential inequality

$$\frac{df}{dt} \le \Delta f + bf,$$ 
on $M\times [0,T]$, where
$\Delta$ refers to a Laplacian at time $t$. Then for any $x\in M$,
$t\in [0,T)$,

$$|f(x,t)| \le c\frac{1}{\sqrt{V}}e^{cHd}\max(1,d^{\frac{n}{2}})(b + l
+ \frac{1}{t})^{\frac{1 + n/2}{2}} e^{cbt}||f_0||_{L^2},$$
where $c$ is a positive constant depending only on $n$ and $d =
\max_{0 \le t \le T} \diam(M,g(t))$, $H = \max_{0 \le t \le T}
\sqrt{||\ric||_{C^0}}$, $f_0 = f(\cdot, 0)$, $V = \min_{0 \le t \le T}
\vol_{g(t)}(M)$.
\end{lemma}

The following remark will give us a condition that will imply
obtaining the Einstein metrics in the limit.

\begin{remark}
\label{remark-remark_speculation}
If $g(t)$ is a solution to $(g_{ij})_t = -2R_{ij} +
\frac{1}{\tau}g_{ij}$, for $t\in [0,\infty)$ such that 
\begin{enumerate}

\item
A curvature operator and a diameter are uniformly bounded along
the flow.
\item
$0 \le R(x,t) \le \frac{n}{2\tau}$ for all $x\in M$ and all
$t\in[0,\infty)$.
\end{enumerate}
Then all the solitons that arise as limits of the subsequences of our
flow $g(t)$ are Einstein metrics with scalar curvatures $R =
\frac{n}{2\tau}$ and $T_{ij}(t)$ converge to zero, uniformly on $M$ as
$t\to\infty$. $T_{ij} = R_{ij} - \frac{R}{n}g_{ij}$ is a traceless
part of the Ricci curvature.
\end{remark}

\begin{proof}[Proof of the Remark]
Notice that now we do not make an assumption that one of the metrics
that we get in a limit is an Einstein metric. Look at the evolution
equation for $r(t) =
\frac{1}{\vol_t(M)}\int_M R dV_t$,
$$\frac{d}{dt}r(t) = \frac{1}{\vol_t(M)}(2\int_M |T|^2 + (1 -
\frac{2}{n})\int_M R(\frac{n}{2\tau} - R) + r(r - \frac{n}{2\tau}).$$
$R \le \frac{n}{2\tau}$ implies $r(t) \le \frac{n}{2\tau}$ and
therefore

\begin{equation}
\label{equation-equation_evolution_r}
\frac{d}{dt}r(t) \ge \frac{2}{\vol_t(M)}\int_M |T|^2 + r(r -
\frac{n}{2\tau}).
\end{equation}
We have proved that in the case of flow $g(t)$, a volume noncollapsing
condition holds for all times $t\ge 0$. $\frac{d}{dt}\ln(\vol_t(M)) =
\frac{n}{2\tau} - r$ and $C_1 \le
\vol_t(M) \le C_2$ give that $\int_0^{\infty} (\frac{n}{2\tau} -
r(t))dt < \infty$. We can integrate the inequality
(\ref{equation-equation_evolution_r}) in $t\in[0,\infty)$. This,
together with the uniform estimates on $\vol_t(M)$ and $r(t)$ give
that
\begin{equation}
\label{equation-equation_estimate_T}
\int_0^{\infty}\int_M |T|^2 dV_t \le C.
\end{equation}
Following the calculations in Hamilton's paper \cite{hamilton1982},
Rugang computed the evolution equation for $T$ under a normalized
Ricci flow (\cite{rugang1993}). In the case of flow
(\ref{equation-equation_our_flow}) we have
\begin{equation}
\label{equation-equation_evolution_T}
\frac{d}{dt}|T|^2 = \Delta|T|^2 - 2|\nabla T|^2 + 4\rem(T)\cdot T +
\frac{4}{n}(R - \frac{n}{2\tau})|T|^2.
\end{equation}
Since the curvature operators of $g(t)$ are uniformly bounded, we
derive from equation (\ref{equation-equation_evolution_T}) that
$$\frac{d}{dt}|T| \le \Delta|T| + C|T|.$$ 
Applying Lemma \ref{lemma-lemma_Moser} to this differential 
inequality and intervals $[t-1,t+1]$ where $t > 1$, we derive
$$|T|^2(x,t) \le ||T||^2(t)_{C^0(M)} \le C(\int_{M_{t-1}} |T|^2),$$
where $M_t = (M,g(t))$. Integrate this inequality in $t\in [k,k+1]$,
for all $k\ge k_0$ and sum up all the inequalities that we get 
this way. We get
\begin{eqnarray*}
\int_{k_0}^{\infty}||T||^2 dt &\le& C\sum_{k\ge k_0}
\int_k^{k+1}(\int_{M_{t-1}} |T|^2) dt
\end{eqnarray*}
\begin{equation}
\label{equation-equation_need_1}
\int_{k_0}^{\infty} ||T||^2 dt \le C\int_{k_0}^{\infty}\int_M |T|^2
dV_{t-1} dt,
\end{equation} 
where $dV_{t-1}$ is a volume form for metric $g(t -
1)$.  $\int_M |T|^2 dV_{t-1} \le C\int_M |T|^2 dV_t$, because
$\frac{d}{dt}\ln\vol_t = \frac{n}{2\tau} - R$ and the 
curvatures of $g(t)$ are uniformly bounded. The right hand side of
inequality (\ref{equation-equation_need_1}) is bounded by a uniform 
constant, because of the estimate (\ref{equation-equation_estimate_T}).  
Therefore $\int_{k_0}^{\infty} ||T||^2dV_t \le C$.

If there exists $(p,t_0)$ such that $|T|^2(p,t_0) > \epsilon$, then
there is a small neighbourhood of $(p,t_0)$ in $M\times [0,\infty)$,
say $U_{\delta}(p,t_0) = B_p(\delta,t_0)\times [t_0 - \delta, t_0 +
\delta]$ such that $|T|^2(x,t) \ge \frac{\epsilon}{2}$ for all
$(x,t)\in U_{\delta}(p,t_0)$. This follows from the fact that in the
case of a Ricci flow, a bound $|\rem| \le C$ implies $|D^kD_t^l \rem|
\le C(k,l)$. Costant $\delta$ does not depend on a point $(p,t_0)\in
M\times [0,\infty)$, since all our bounds and estimates are uniform.

If there existed $\epsilon > 0$ and a sequence of points $(p_i,t_i)\in
M\times [0,\infty)$, with $t_i\to\infty$ such that $|T(p_i,t_i)| \ge
\epsilon$ then we would have that $||T||_{C^0} \ge \frac{\epsilon}{2}$
for all $t\in [t_i - \delta, t_i + \delta]$ and for all $i$. This
would imply $C \ge \int_0^{\infty}||T||^2 dV_t \ge \sum_{i=0}^{\infty}
\epsilon \delta = \infty$. This is impossible. Therefore,
$||T||_{C^0(M_t)}\to 0$ as $t\to\infty$.

$\frac{d}{dt}\ln(\vol_t) = \frac{n}{2\tau} - R \ge 0$ for all $t$
imply that there exists a finite $\lim_{t\to\infty}\vol_t$ for every
$x\in M$ (otherwise we can argue as in the previous paragraph). If we
integrate this equation in $t\in [0,\infty)$, we will get that
$\int_0^{\infty} (\frac{n}{2\tau} - R) dt < \infty$. As in the case
for a traceless part of the Ricci curvature $T$, we can conclude that
$\lim_{t\to\infty}R = \frac{n}{2\tau}$ uniformly on $M$.

We can conclude that under the assumptions given at the beginning of
this remark, for every sequence $t_i\to\infty$ we can find a
subsequence such that $g(t_i + t) \to h(t)$, where $h(t)$ is an Einstein
soliton with scalar curvature $\frac{n}{2\tau}$. We also know that
$R_{ij} - \frac{1}{2\tau}g_{ij} \to 0$ as $t\to\infty$, uniformly on
$M$ and that there exists $\lim_{t\to\infty}\vol_t$. 
\end{proof}

To conclude, we have proved a sequential convergence of a solution of
a $\tau$-flow towards solitons (generalizations of Einstein metrics),
under uniform curvature and diameter assumptions. We still do not know
whether we get a unique soliton (up to diffeomorphisms) in the limit
or not. All observations in this subsection are in favour of the
uniqueness of a soliton in the limit. 

\end{subsection}

\end{section}


\begin{thebibliography}{20}

\bibitem{cao1985} H.D.Cao: Deformation of Kahler metrics to
Kahler-Einstein metrics on comapct Kahler manifolds; Invent. math. 81
(1985) 359--372.

\bibitem{colding} J. Cheeger, T. Colding: On the structure of spaces
with Ricci curvature bounded below.

\bibitem{colding1997} T. Colding: Ricci curvature and volume
convergence; Annals of Mathematics, 145 (1997), 477--501.

\bibitem{deturck1983} D. Deturck: Deforming metrics in the direction
of their ricci tensors; J. Diff. Geom. 18 (1983), 157--162.

\bibitem{glickenstein2002} D. Glickenstein: Precompactness of
solutions to the Ricci flow in the absence of injectivity radius
estimates; preprint arXiv:math.DG/0211191 v2.

\bibitem{hamilton1982} R. Hamilton: Three-manifolds with positive
Ricci curvature, Journal of Differential Geometry 17 (1982) 225--306.

\bibitem{hamilton1995} A compactness property for solutions of the
Ricci flow, Amer.J.Math. 117 (1995) 545--572.

\bibitem{hamiltonMA} The formation of singularities in the Ricci flow,
Surveys in Differential Geometry, vol. 2, International Press,
Cambridge, MA (1995) 7--136.

\bibitem{hamilton1999} R. Hamilton:Non-singular solutions of the Ricci
flow on 3 manifolds, Communications in analysis and geometry vol. 7
(1999) 695--729.

\bibitem{perelman2002} G. Perelman: The entropy formula for the Ricci
flow and its geometric applications, preprint.

\bibitem{rothaus1981} Rothaus: Logarithmic Sobolev inequality and
spectrum, Journal of Dif. Anal. 42 (1981) 109--120.

\bibitem{rugang1993} Rugang Ye: Ricci flow, Einstein metrics and space
forms, Transactions of the american mathematical society, volume 338,
number 2 (1993) 871--895.

\end{thebibliography}
\end{document}